\documentclass[12pt]{article}

\newtheorem{theorem}{Theorem}
\newtheorem{corollary}{Corollary}

\newtheorem{remark}{Remark}
\newtheorem{example}{Example}

\newcommand{\mymod}[3]{$#1 \equiv #2 \,\, (\!\!\!\!\mod #3)$}

\newcommand{\JK}[2]{
\begin{array}{lll}
J & = & \{ #1 \} \\
K & = & \{ #2 \} \\
\end{array}
}

\usepackage{epsfig}
\usepackage{url}
\usepackage{amssymb}
\usepackage{amsmath}
\usepackage{amsfonts}

\def\qed{\begin{flushright} $\Box$ \end{flushright}}

\def\Dbar{\leavevmode\lower.6ex\hbox to 0pt{\hskip-.23ex \accent"16\hss}D}

\begin{document}

{\bf\LARGE
\begin{center}
New results on D-optimal Matrices\footnote{This research is
supported by NSERC grants.}
\end{center}
}

{\Large
\begin{center}
Dragomir {\v{Z}}. {\Dbar}okovi{\'c}\footnote{University of Waterloo,
Department of Pure Mathematics, Waterloo, Ontario, N2L 3G1, Canada
e-mail: \url{djokovic@math.uwaterloo.ca}}, Ilias S.
Kotsireas\footnote{Wilfrid Laurier University, Department of Physics
\& Computer Science, Waterloo, Ontario, N2L 3C5, Canada, e-mail:
\url{ikotsire@wlu.ca}}
\end{center}
}

\begin{abstract}
We construct a number of new $(v;r,s;\lambda)$ supplementary
difference sets (SDS) with $v$ odd and $\lambda = (r+s)-(v-1)/2$. In particular, these
give rise to D-optimal matrices of the four new orders $206,
242, 262, 482$, constructed here for the first time.
\end{abstract}

\section{Introduction}

\noindent Let $v$ be an odd positive integer and consider
$(-1,1)$-matrices $H$ of order $2v$. Ehlich's bound states that
$$
    \det(H) \leq 2^v (2v-1)(v-1)^{v-1}.
$$
D-optimal matrices are $2v \times 2v$ $(-1,1)$-matrices that attain
Ehlich's bound, i.e. they have maximal determinant. Ehlich also
proved that if $A$ and $B$ are circulant $(-1,1)$-matrices of order
$v$ such that
\begin{equation}
    A A^T + B B^T = 2(v-1) I_v + 2J_v
\label{D-optimal-matrix-equation}
\end{equation}
(where $I_v$ is the $v \times v$ identity matrix and $J_v$ is the $v
\times v$ matrix with all elements equal to $+1$), then the matrix
\begin{equation}
H = \left(
\begin{array}{cc}
A & B \\
-B^T & A^T
\end{array}
\right)
\label{D-optimal-from-two-circulants}
\end{equation}
has maximal determinant. Such $A$ and $B$ can be constructed by
using cyclic SDSs, say $(X,Y)$,
with parameters $(v;r,s;\lambda)$. Thus $|X|=r$, $|Y|=s$ and
$\lambda =r+s-\frac{v-1}{2}$. We say that D-optimal matrices
(\ref{D-optimal-from-two-circulants}) are of {\it circulant type}
and that such SDSs are {\it D-optimal}.
By pre- and post-multiplying equation
(\ref{D-optimal-matrix-equation}) with $J_v$, one obtains that
$a^2+b^2=4v-2$, where $a$ and $b$ are  row sums of $A$ and $B$,
respectively. Each integer solution $(a,b)$ of this Diophantine
equation can be used to obtain feasible parameters $r$ and $s$
for the required SDS. Namely, we need $a=v-2r$ and $b=v-2s$.
By a normalization, we may assume that $0<a\le b$ which implies
that $r\ge s$.

\noindent A comprehensive table of all odd $v < 100$ for which
D-optimal SDSs are known, can be found in \cite{Djokovic:1997}.
Additional values of $v < 100$ for which D-optimal matrices of order
$2v$ have been found subsequently, can be found in
\cite{KharaghaniOrrick:2007}. There are two infinite series of
D-optimal matrices, one for $v=q^2+q+1$ where $q$ is a prime power,
see \cite{KKS:1991:Infinite-Series}, and one for $v=2q^2+2q+1$ where
$q$ is an odd prime power, see \cite{Whiteman:1990:Infinite-Series}.
Apart from these two infinite series, the only odd values of $v >
100$ for which D-optimal matrices of order $2v$ are currently known
are $v = 113$, see \cite{Djokovic:1997,Gysin:1997}, $v =
145$, see \cite{Djokovic:1997,Gysin:Thesis:1997}, $v = 157$,
see \cite{Gysin:Thesis:1997,GysinSeberry:1998} and $v = 181$, see
\cite{Gysin:Thesis:1997,KSWX:JSPI:1997}. We point out that the cases
$v = 157$ and $v = 181$ were overlooked in \cite{KharaghaniOrrick:2007}.

\noindent In this paper we construct several new D-optimal SDSs for
seven new parameter sets:
$$
\begin{array}{lll}
(63;29,24;22)  & (93;45,37;36)  & (103;48,42;39) \\
(103;46,43;38) & (121;55,51;46) & (131;61,55;51) \\
(241;120,105;105). & &
\end{array}
$$
In particular, this means that we have constructed the first
examples of D-optimal matrices of orders $206$, $242$, $262$ and
$482$. D-optimal matrices of order $126$ and $186$ constructed
previously in \cite{Djokovic:1991} and \cite{Djokovic:1997}
were also derived from D-optimal SDS but with
different parameters, namely $(63;27,25;21)$ and $(93;42,38;34)$,
respectively. However the
smallest order $n \equiv 2 \pmod{4}$ for which a D-optimal matrix is
still unknown remains $n = 138$. For a comprehensive list of unknown
D-optimal SDSs see Section 5.

\noindent We also prove two theoretical results for D-optimal
matrices. The first is a generalization of the horizontal and
vertical constraints for D-optimal matrices, recently proved in
\cite{KotsireasPardalos:2011}. The second is the property that if a
D-optimal SDS solution is formed by taking union of cosets of a
certain subgroup of the group of units $Z_v^\star$, then the power
spectral density of the corresponding sequences must be constant on
the cosets. (See below the definition of the power
spectral density.) This property can be restated as saying that
several power spectral density values corresponding to the first
rows of the circulant matrices $A$, $B$ in
(\ref{D-optimal-from-two-circulants}), must be equal. Both our
theoretical results can potentially be used as subroutines to
discard candidate sequences or candidates SDSs, within algorithmic
schemes that search for D-optimal matrices.

\noindent Let us recall the definition of the power spectral
density (PSD) of a complex sequence  $X = x_0,x_1,\ldots,x_{v-1}$.
Let $\omega = e^{\frac{2\pi i}{v}}$ be a primitive $v$-th
root of unity and define the discrete Fourier transform (DFT) of
$X$ by
$$
    DFT_X(k) = \sum_{\alpha = 0}^{v-1} x_{\alpha} \omega^{\alpha k},
    \quad k = 0,\ldots,v-1.
$$
Then the power spectral density values $PSD_X(k)$ are defined by
$$
    PSD_X(k) = |DFT_X(k)|^2, \quad k = 0,1,\ldots,v-1.
$$

\noindent From now on we work only with $\{\pm1\}$-sequences.
The well-known symmetry properties
$$
 DFT_X(v-k) = \overline{DFT_X(k)}  \mbox{ and }
 PSD_X(v-k) = PSD_X(k), \quad k = 0,1,\ldots,v-1
$$
are used within most algorithmic schemes that search for D-optimal
matrices.

\section{Generalized horizontal and vertical constraint}

\noindent We state and prove a generalized constraint, that yields
the horizontal and vertical constraints proved in
\cite{KotsireasPardalos:2011} for \mymod{v}{0}{3}, as
corollaries. In particular we show that actually any divisor $d$ of
$v$ gives rise to a Diophantine constraint.

\begin{theorem}
Let the first rows of the circulant matrices $A,B$ in
(\ref{D-optimal-from-two-circulants}) be denoted as
$[a_0,\ldots,a_{v-1}]$ and $[b_0,\ldots,b_{v-1}]$ respectively,
where $a_i,b_i \in \{-1,+1\}, i=0,\ldots,v-1$. Let $v = dm$ and set
$$
\begin{array}{c}
    A_j = a_j + a_{j+d} + \cdots + a_{j+(m-1)d} \\
    B_j = b_j + b_{j+d} + \cdots + b_{j+(m-1)d} \\
\end{array}
$$
for $j=0,\ldots,d-1$. If $A,B$ define a circulant D-optimal matrix,
then
\begin{equation}
    \sum_{j=0}^{d-1} ( A_j^2 + B_j^2 ) = 2(v+m-1)
\label{Equation_General_Constraint}
\end{equation}
and
\begin{equation}
    \sum_{k<l} ( A_k A_l + B_k B_l ) = v-m.
    \label{Zbir-AA-BB}
\end{equation}

\label{Theorem:Generalized_Vertical_Constraint}
\end{theorem}
\noindent {\bf Proof} \\
By abuse of notation we denote by $A,B$ the circulant matrices and
the sequences that are specified by their first rows. An argument
similar to the one given in \cite{FGS:2001}, shows that the hypothesis
that $A,B$ define a circulant D-optimal matrix implies that the sum
$PSD_A(s)+PSD_B(s)$, $s\ne0$, is equal to the constant $2v-2$, i.e.
$$
    \left( \sum_{k=0}^{v-1} a_k \cos\frac{2\pi ks}{v} \right)^2 +
    \left( \sum_{k=0}^{v-1} a_k \sin\frac{2\pi ks}{v} \right)^2 +
    \left( \sum_{k=0}^{v-1} b_k \cos\frac{2\pi ks}{v} \right)^2 +
    \left( \sum_{k=0}^{v-1} b_k \sin\frac{2\pi ks}{v} \right)^2 =
    2v - 2
$$
for all $s \in \{ 1,2,\ldots,v-1 \}$. (For $s=0$ the sum is
equal to $4v-2$.) Considering the particular
values $s=mr$ for $r=1,2,\ldots,d-1$ we obtain:
\begin{equation}
    \left( \sum_{k=0}^{v-1} a_k \cos\frac{2\pi kr}{d} \right)^2 +
    \left( \sum_{k=0}^{v-1} a_k \sin\frac{2\pi kr}{d} \right)^2 +
    \left( \sum_{k=0}^{v-1} b_k \cos\frac{2\pi kr}{d} \right)^2 +
    \left( \sum_{k=0}^{v-1} b_k \sin\frac{2\pi kr}{d} \right)^2 =
    2v - 2.
    \label{PSD_s_eq_m}
\end{equation}
Noting that
$$
\begin{array}{ccl}
    \left( \displaystyle\sum_{k=0}^{v-1} a_k \cos\frac{2\pi kr}{d} \right)^2 +
    \left( \displaystyle\sum_{k=0}^{v-1} a_k \sin\frac{2\pi kr}{d} \right)^2
    & = & \left( \displaystyle\sum_{j=0}^{d-1} A_j \cos\frac{2\pi jr}{d} \right)^2 +
    \left( \displaystyle\sum_{j=0}^{d-1} A_j \sin\frac{2\pi jr}{d} \right)^2 \\
    & = & \displaystyle\sum_{j=0}^{d-1} A_j^2 + 2 \sum_{k<l} A_k A_l \cos\frac{2\pi
(l-k)r}{d}
\end{array}
$$
(and similarly for the B sequence), we see that (\ref{PSD_s_eq_m})
can be rewritten as
\begin{equation}
    \sum_{j=0}^{d-1} ( A_j^2 + B_j^2 ) + 2 \sum_{k<l} ( A_k A_l + B_k B_l ) \cos\frac{2\pi(l-k)r}{d} = 2v-2.
    \label{Zbir-A-B-r}
\end{equation}
The hypothesis that $A,B$ define a circulant D-optimal matrix
implies that the sum of the squares of the sums of the elements of
the sequences $A$ and $B$ is equal to the constant $4v-2$, i.e.
$$
    \left( \displaystyle\sum_{k=0}^{d-1} A_k \right)^2 +
    \left( \displaystyle\sum_{k=0}^{d-1} B_k \right)^2
    = 4v-2,
$$
which can be rewritten as
\begin{equation}
    \sum_{j=0}^{d-1} ( A_j^2 + B_j^2 ) +
    2 \sum_{k<l} ( A_k A_l + B_k B_l ) = 4v-2.
    \label{Zbir-A-B}
\end{equation}
By taking the sum of (\ref{Zbir-A-B-r}) over $r=1,\ldots,d-1$
and of (\ref{Zbir-A-B}), we obtain the equality
(\ref{Equation_General_Constraint}).
Now (\ref{Zbir-AA-BB}) follows from (\ref{Zbir-A-B}).
\qed

\noindent Taking the difference of (\ref{Zbir-A-B}) and 
(\ref{Zbir-A-B-r}), we obtain
$$
    \sum_{0 \leq k < l < d} (A_k A_l + B_k B_l) \left( 1 -
\cos\frac{2\pi(l-k)r}{d} \right) = v,
$$
which can be rewritten as
\begin{equation}
    \sum_{0 \leq k < l < d} (A_k A_l + B_k B_l) \sin^2\frac{\pi(l-k)r}{d}
    = \displaystyle\frac{v}{2}.
\label{Equation_Generalized_Vertical_Constraint}
\end{equation}

\noindent We remark that assuming \mymod{v}{0}{3} and
applying theorem \ref{Theorem:Generalized_Vertical_Constraint} for
$d=3$ and $d=\frac{v}{3}$ yields immediately the vertical and horizontal constraints proved in \cite{KotsireasPardalos:2011}.

\begin{corollary} $ $ \\
\noindent (a) Assuming \mymod{v}{0}{3} and taking
$d=3$ in theorem \ref{Theorem:Generalized_Vertical_Constraint} we
obtain the vertical constraint:
$$
        A_0^2 + A_1^2 + A_2^2 + B_0^2 + B_1^2 + B_2^2  = 8 m - 2
$$
where $A_j = \displaystyle\sum_{i=0}^{m-1} a_{3i+j}, B_j =
\displaystyle\sum_{i=0}^{m-1} b_{3i+j}, j = 0,1,2$.

\noindent (b) Assuming \mymod{v}{0}{3} and taking
$d=\frac{v}{3}$ in theorem \ref{Theorem:Generalized_Vertical_Constraint} we
obtain the horizontal constraint:
$$
     \sum_{i=0}^{d-1} ( A_i^2 + B_i^2 ) = 2v+4
$$
where $A_i = a_i+a_{i+d}+a_{i+2d}, B_i = b_i+b_{i+d}+b_{i+2d},
i=0,\ldots,d-1$.
\end{corollary}

\noindent Let us also illustrate the interesting phenomenon that
in case the irrationalities (that the presence of the sine function
in (\ref{Equation_Generalized_Vertical_Constraint}) potentially
entails) do not cancel out, then one obtains more than one
Diophantine equation.

\begin{corollary} Assuming \mymod{v}{0}{5} and taking
$d=5$ and $r=1$ in 
(\ref{Equation_Generalized_Vertical_Constraint}) 
we obtain the constraints:
$$
        \sum_{0 \leq k < l < 5} (A_k A_l + B_k B_l) = 4m
        \mbox{ and }
        \sum_{l-k \, \in \{2,3\}} (A_k A_l + B_k B_l)
        -
        \sum_{l-k \, \in \{1,4\}} (A_k A_l + B_k B_l) = 0
$$
where $A_j = \displaystyle\sum_{i=0}^{m-1} a_{5i+j}, B_j =
\displaystyle\sum_{i=0}^{m-1} b_{5i+j}, j = 0,1,\ldots,4$.

\end{corollary}
\noindent {\bf Proof} \\
For $d=5, v=5m$, the summation in
(\ref{Equation_Generalized_Vertical_Constraint}) for $r=1$
extends over the $10$ possible $(k,l)$ pairs with
$0 \leq k < l \leq 4$ and in each case we have that
$\sin^2\displaystyle\frac{\pi (l-k)}{5}$ is either
equal to $\sin^2\displaystyle\frac{\pi}{5} =
\displaystyle\frac{5-\sqrt{5}}{8}$ or equal to
$\sin^2\displaystyle\frac{2\pi}{5} =
\displaystyle\frac{5+\sqrt{5}}{8}$. Since $1$ and $\sqrt{5}$ are
linearly independent over the rationals, we obtain the aforementioned
constraints. \qed

\section{Power spectral density constancy over orbits}

\noindent First we introduce some notations and terminology. Let
$Z_v$ be the ring of integers $\mod v$, i.e $Z_v= \{
0,1,\ldots,v-1\}$. Let $Z_v^\star$ be the group of invertible
elements of $Z_v$, i.e. $Z_v^\star = \{ k \in Z_v : \gcd(k,v) = 1
\}$. The order of $Z_v^\star$ is equal to $\phi(v)$, where $\phi$
denotes the Euler totient function. Let $H \leqslant Z_v^\star$ be a
subgroup of $Z_v^\star$. Then $H$ acts on $Z_v$ and we denote the
orbits of this action by
$$
    {\cal O}_1 = \{ 0 \}, {\cal O}_2, \ldots {\cal O}_m.
$$
Thus we have the disjoint union relationship $Z_v = {\cal O}_1 \cup
{\cal O}_2 \cup \cdots \cup {\cal O}_m$.

\noindent The method used in this paper as well as in
\cite{Djokovic:1991,Djokovic:1997,Gysin:Thesis:1997,Gysin:1997,GysinSeberry:1998,KSWX:JSPI:1997}
constructs solutions for circulant type D-optimal matrices by expressing the
corresponding SDSs as unions of certain orbits associated to a
suitable subgroup of $Z_v^\star$. This special structure of these
solutions implies certain constraints on the possible range of
values of the power spectral densities of the two sequences
associated to the SDS. Specifically, we state and prove the fact
that the power spectral densities remain constant over the orbits.

\noindent Let $(X,Y)$ be an SDS of $Z_v$ with parameters
$(v;r,s;\lambda)$, with $v$ odd and $\lambda = r+s-\frac{v-1}{2}$,
corresponding to a circulant D-optimal matrix. In particular, we
must have $|X|=r$ and $|Y|=s$. Assume that
\begin{equation}
    X = \displaystyle\bigcup_{j \in J} {\cal O}_j, \quad
    Y = \displaystyle\bigcup_{k \in K} {\cal O}_k
\label{X_and_Y_as_Unions_of_orbits}
\end{equation}
for some subsets $J,K$ of $\{ 1,2,\ldots,m \}$.

\noindent By abuse of notation, let $X$ also denote the sequence
$x_0,x_1,\ldots,x_{v-1}$ where
$$
    x_i =
    \left\{
    \begin{array}{rcc}
    1 & \mbox{if} & i \not\in X \\
    -1 & \mbox{if} & i \in X \\
    \end{array}
    \right.
$$
and define similarly the sequence $Y = y_0,y_1,\ldots,y_{v-1}$.

\begin{theorem}
If $k$ and $k'$ belong to the same orbit ${\cal O}_r \subseteq Z_v$
and the sequence $X$ is as defined above, then
$$
    PSD_X(k) = PSD_X(k').
$$
\label{Theorem:PSD_Equal_on_Orbits}
\end{theorem}
\noindent {\bf Proof} \\
Since $k,k'$ belong to the same orbit ${\cal O}_r$ we have
$k'=ks$ for some $s \in H$. Since the set $X$ is a union of orbits,
we also have $x_{\alpha s}=x_\alpha$ for all $\alpha\in Z_v$ and
$s\in H$. Hence we obtain
$$
    PSD_X(k') = |DFT_X(k')|^2 =
    \left| \sum_{\alpha = 0}^{v-1} x_{\alpha} \omega^{\alpha k'} \right|^2 =
    \left( \sum_{\alpha = 0}^{v-1} x_{\alpha} \omega^{\alpha k'} \right)
    \cdot
    \left( \sum_{\beta = 0}^{v-1} x_{\beta} \omega^{-\beta k'} \right)
$$
because $x_{\alpha} \in \{-1,+1\}$ and $\bar{\omega} = \omega^{-1}$.
By expanding we have
$$
\begin{array}{ccc}
    PSD_X(k') & = &
    \displaystyle\sum_{\alpha = 0, \beta =0 }^{v-1} x_{\alpha} x_{\beta} \, \omega^{(\alpha - \beta)k'} =
    \displaystyle\sum_{\alpha = 0, \beta =0 }^{v-1} x_{\alpha} x_{\beta} \, \omega^{(\alpha - \beta)ks} =
    \displaystyle\sum_{\alpha = 0, \beta =0 }^{v-1} x_{\alpha s} x_{\beta s} \, \omega^{(\alpha - \beta)ks}
= \\
       & = & \displaystyle\sum_{\alpha' = 0, \beta' =0 }^{v-1} x_{\alpha'} x_{\beta'} \, \omega^{(\alpha' - \beta')k} =
    \left| \displaystyle\sum_{\alpha' = 0}^{v-1} x_{\alpha'} \omega^{\alpha' k} \right|^2 =
    |DFT_X(k)|^2 = PSD_X(k). \\
\end{array}
$$
Here we have used the fact that when $\alpha$ ranges through $Z_v$,
so does $\alpha' = \alpha s$, when $s \in H$.

\qed

\begin{remark}
\noindent Theorem \ref{Theorem:PSD_Equal_on_Orbits} says that the
PSD is constant on the orbits of $H$. In fact it can easily be shown
that PSD is constant on the orbits of the possibly larger group
$H^\star = H \cup (-1)H$.

\noindent Theorem \ref{Theorem:PSD_Equal_on_Orbits} has some
important computational consequences when one is using the power
spectral density property as a preprocessor to discard candidate
sequences and SDSs from a search for D-optimal matrices. In
particular, it is clear that one needs to compute only one PSD value
from each orbit ${\cal O}_j$ for $2\le j\le m$. 
Therefore, one can reduce the
number of PSD values calculated per SDS from $\frac{v-1}{2}$ to
at most $m-1$.

\end{remark}

\begin{example}
We illustrate theorem \ref{Theorem:PSD_Equal_on_Orbits} with a
specific solution for $v = 131$, a prime. The automorphism group $G$
of the cyclic group $Z_{131}$ is cyclic of order $\phi(131)=130$.
The subgroup $H$ of $G$ that was used was one of order $5$. Its
generator is $53$, $H = \{ 1,53,58,61,89 \}$. 
The $H$-orbits in $Z_v$ are $\{0\}$ and the 26 cosets of $H$ in $G$, 
thus $m=27$. We enlarge this group
by using the element $-1$ ($=130$ in this case). Then we get a
subgroup $H^\star = \{ 1,42,53,58,61,70,73,78,89,130 \}$ of order
$10$. There are $13$ cosets of $H^\star$ in $G$, all of size $10$ as
well. Theorem \ref{Theorem:PSD_Equal_on_Orbits} implies that
$$
    PSD_X(k') = PSD_X(k) \mbox{ and } PSD_Y(k') = PSD_Y(k), \quad \forall k' \in H^\star \cdot k.
$$
In this example, one can compute $\frac{m-1}{2} = 13$ 
PSD values for each sequence instead of
$\frac{v-1}{2}=65$ values for each sequence.
\end{example}

\section{The new circulant D-optimal matrices and SDSs}

\noindent We give D-optimal SDSs for the seven new parameter sets
listed in the introduction. For some of them we give two or three
non-equivalent SDSs. Non-equivalence of SDSs was established in
all cases by an implementation of the method described in \cite{Djokovic:2011}. \\

\noindent We define the notation for the orbits of the action of the
subgroup that we use to construct the solutions below. Consider a
fixed subgroup $H$ of order $h$ of $Z_v^\star$. Clearly $h$ must
divide $|Z_v^\star| = \phi(v)$. Denote by $H \cdot k$ the orbit of
$H$ in $Z_v$ through the point $k$, where $\cdot$ is multiplication
$\mod \,\, v$. We refer to the orbit $H \cdot 0 = \{ 0 \}$ as the
{\it trivial orbit}. The orbit $H \cdot 1$ is just the subgroup $H$
itself. In general the size of an orbit may be any divisor of $|H|$
and if $v$ is a prime then every nonzero orbit is just a coset of
$H$ in $Z_v^\star$ and so will have size $|H|$.

\noindent The notation
\begin{equation}
      X = \bigcup_{j \in J} H \cdot j
      \mbox{  and  }
      Y = \bigcup_{k \in K} H \cdot k
        \label{Notation:X_Y_J_K_SDS}
\end{equation}
will be used below to present all the solutions found, in fact each
solution will be given only via the two indexing sets $J$ and $K$.
When $X$, $Y$ are defined as above, then $(X,Y)$ will be an
SDS$(v;r,s;\lambda)$ with $r = \sum_{j \in J} |H \cdot j|$, $s =
\sum_{k \in K} |H \cdot k|$ and $\lambda = r+s-\frac{v-1}{2}$.

\subsection{$v = 63$}

In this case the subgroup $H$ is trivial and we write the blocks
$X$ and $Y$ explicitly:
\begin{eqnarray*}
X &=& \{0,1,2,3,4,6,7,11,12,13,14,20,21,22,25,26,27,30,33,35,36,38,39,42,46,48,50,53,57\}, \\
Y &=& \{0,1,3,5,7,8,10,11,13,14,16,18,19,23,30,33,34,35,39,40,48,52,54,56\}.
\end{eqnarray*}
This is an SDS with parameters $(63;29,24;22)$ which gives rise
to a D-optimal matrix of order $2 \times 63 = 126$.

\subsection{$v = 93$}

Consider the subgroup $H = \{ 1,25,67 \}$ of order $3$, of
$Z_{93}^\star$. For the convenience of the reader we give below the
$33$ orbits of the action of $H$ on $Z_{93}$.
$$
\begin{array}{lll}
H \cdot 0 = \{0\}         & H \cdot 1 = \{1,25,67\}   & H \cdot 2 = \{2,41,50\}  \\
H \cdot 3 = \{3,15,75\}   & H \cdot 4 = \{4, 7, 82\}  & H \cdot 5 = \{5,32,56\}  \\
H \cdot 6 = \{6,30,57\}   & H \cdot 8 = \{8,14,71\}   & H \cdot 9 = \{9,39,45\}  \\
H \cdot 10 = \{10,19,64\} & H \cdot 11 = \{11,86,89\} & H \cdot 12 = \{12,21,60\}  \\
H \cdot 13 = \{13,34,46\} & H \cdot 16 = \{16,28,49\} & H \cdot 17 = \{17,23,53\}  \\
H \cdot 18 = \{18,78,90\} & H \cdot 20 = \{20,35,38\} & H \cdot 22 = \{22,79,85\}  \\
H \cdot 24 = \{24,27,42\} & H \cdot 26 = \{26,68,92\} & H \cdot 29 = \{29,74,83\}  \\
H \cdot 31 = \{31\}       & H \cdot 33 = \{33,72,81\} & H \cdot 36 = \{36,63,87\}  \\
H \cdot 37 = \{37,61,88\} & H \cdot 40 = \{40,70,76\} & H \cdot 43 = \{43,52,91\}  \\
H \cdot 44 = \{44,65,77\} & H \cdot 47 = \{47,59,80\} & H \cdot 48 = \{48,54,84\}  \\
H \cdot 51 = \{51,66,69\} & H \cdot 55 = \{55,58,73\} & H \cdot 62 = \{62\}  \\
\end{array}
$$
We give three non-equivalent SDSs with parameters $(93;45,37;36)$,
via their index sets $J,K$ to be used in
(\ref{Notation:X_Y_J_K_SDS}), which give rise to D-optimal matrices
of order $2 \times 93 = 186$.
$$
\begin{array}{l}
\JK{ 2,5,8,9,10,12,13,24,33,36,37,40,43,47,55  }{ 3,4,5,6,16,22,24,26,33,36,40,44,62   } \\
 \\
\JK{ 3,6,8,10,11,12,13,17,18,22,29,33,37,43,55 }{ 3,4,5,12,13,16,17,26,36,37,40,51,62  } \\
 \\
\JK{ 3,5,8,10,11,12,13,16,17,20,24,29,33,48,55 }{ 5,9,11,12,13,18,22,24,29,40,43,51,62 } \\
\end{array}
$$

\subsection{$v = 103$}

Consider the subgroup $H = \{ 1,46,56 \}$ of order $3$, of
$Z_{103}^\star$. For the convenience of the reader we give below the
$35$ orbits of the action of $H$ on $Z_{103}$.
$$
\begin{array}{lll}
H \cdot 0 = \{ 0 \}        & H \cdot 1 = \{1,46,56\}  & H \cdot 2 = \{2, 9, 92\}     \\
H \cdot 3 = \{3,35,65\}    & H \cdot 4 = \{4,18,81\}  & H \cdot 5 = \{5,24,74\}      \\
H \cdot 6 = \{6,27,70\}    & H \cdot 7 = \{7,13,83\}    & H \cdot 8 = \{8,36,59\}    \\
H \cdot 10 = \{10,45,48\}  & H \cdot 11 = \{11,94,101\} & H \cdot 12 = \{12,37,54\}  \\
H \cdot 14 = \{14,26,63\}  & H \cdot 15 = \{15,16,72\}  & H \cdot 17 = \{17,25,61\}  \\
H \cdot 19 = \{19,34,50\}  & H \cdot 20 = \{20,90,96\}  & H \cdot 21 = \{21,39,43\}  \\
H \cdot 22 = \{22,85,99\}  & H \cdot 23 = \{23,28,52\}  & H \cdot 29 = \{29,79,98\}  \\
H \cdot 30 = \{30,32,41\}  & H \cdot 31 = \{31,87,88\}  & H \cdot 33 = \{33,76,97\}  \\
H \cdot 38 = \{38,68,100\} & H \cdot 40 = \{40,77,89\}  & H \cdot 42 = \{42,78,86\}  \\
H \cdot 44 = \{44,67,95\}  & H \cdot 47 = \{47,57,102\} & H \cdot 49 = \{49,66,91\}  \\
H \cdot 51 = \{51,75,80\}  & H \cdot 53 = \{53,69,84\}  & H \cdot 55 = \{55,58,93\}  \\
H \cdot 60 = \{60,64,82\}  & H \cdot 62 = \{62,71,73\}  & \\
\end{array}
$$
We give three non-equivalent SDSs with parameters $(103;48,42;39)$, via their index
sets $J,K$ to be used in (\ref{Notation:X_Y_J_K_SDS}), which give
rise to D-optimal matrices of order $2 \times 103 = 206$.
$$
\begin{array}{l}
\JK{2, 4, 8, 10, 12, 14, 17, 19, 30, 33, 42, 44, 47, 51, 60, 62}{1, 2, 3, 14, 20, 21, 30, 33, 38, 40, 42, 44, 53, 60} \\
 \\
\JK{1, 2, 6, 7, 8, 11, 15, 20, 29, 30, 31, 38, 42, 44, 51, 62}{1, 4, 11, 12, 19, 21, 23, 29, 31, 38, 40, 47, 51, 53} \\
 \\
\JK{3, 5, 8, 10, 12, 14, 15, 21, 22, 29, 30, 44, 47, 51, 60, 62}{1, 2, 5, 14, 15, 19, 23, 29, 31, 33, 40, 42, 51, 55}
\end{array}
$$
We give one SDS with parameters $(103;46,43;38)$, via its index sets
$J,K$ to be used in (\ref{Notation:X_Y_J_K_SDS}), which gives rise
to a D-optimal matrix of order $2 \times 103 = 206$.
$$
\begin{array}{l}
\JK{0, 2, 7, 10, 11, 12, 15, 17, 19, 29, 31, 33, 38, 40, 42, 47}{0, 8, 10, 15, 20, 21, 22, 23, 33, 38, 40, 47, 49, 53, 55}
\end{array}
$$

\subsection{$v = 121$}

Consider the subgroup $H = \{ 1,3,9,27,81 \}$ of order $5$, of
$Z_{121}^\star$. We give one SDS with parameters $(121;55,51;46)$
which gives rise to a D-optimal matrix of order $2 \times 121 =
242$.
$$
\begin{array}{l}
\JK{1,2,5,13,16,19,31,34,35,61,76}{0,7,8,10,13,16,22,25,26,40,76}
\end{array}
$$

\subsection{$v = 131$}

Consider the subgroup $H = \{ 1,53,58,61,89 \}$ of order $5$, of
$Z_{131}^\star$. We give two non-equivalent SDSs with parameters
$(131;61,55;51)$, which give rise to D-optimal matrices of order $2
\times 131 = 262$.
$$
\begin{array}{l}
\JK{0,1,12,14,18,22,27,29,33,36,38,42,44}{2,4,8,11,12,17,22,33,36,38,42} \\
 \\
\JK{0,1,3,4,6,12,14,18,21,22,33,36,38}{2,3,11,12,17,18,19,27,29,38,42} \\
\end{array}
$$

\subsection{$v = 241$}

\noindent Consider the subgroup $H = \{
1,15,24,54,87,91,94,98,100,119,160,183,205,225,231 \}$ of order
$15$, of $Z_{241}^\star$. We give one SDS with parameters $(241;120,105;105)$,
which give rises to a D-optimal matrix of order $2 \times 241 = 482$.
$$
\begin{array}{l}
\JK{ 3, 4, 5, 6, 7, 10, 13, 38 }{ 3, 5, 7, 11, 19, 35, 38 }
\end{array}
$$

\subsection{The Algorithm}

\noindent We now furnish a succinct description of the algorithm used in this paper.
First we record the periodic autocorrelation function (PAF) values of the two candidate sequences 
$A = [a_0,\ldots,a_{v-1}]$ and $B = [b_0,\ldots,b_{v-1}]$, keeping in mind that the relationship 
$$
        PSD_A(k) + PSD_B(k) = 2v-2
$$
must hold for $k = 1,\ldots,v-1$. This relationship is referred to as the PSD criterion and its 
importance lies in the fact that if for a certain value of $k$ one of the $PSD$ values is larger 
than $2v-2$, then the corresponding sequence can be discarded, because of the non-negativity of the 
$PSD$ values. The algorithm proceeds by generating randomly several millions of pairs $(A,B)$ of 
$\pm 1$ sequences of length $v$ that satisfy the PSD criterion. The $A$ sequences are stored into a file 
and encoded by the concatenation of their PAF values on $H^\star$-orbits. 
The $B$ sequences are stored into another file and encoded by the concatenation of their PAF values on $H^\star$-orbits, but where each PAF value is subtracted from the constant $\lambda$.
Then a distributed sorting algorithm is used to sort the two files in ascending order. Subsequently a linear time
algorithm is used to detect any potential matches between the two sorted files. When a match is detected, another
algorithm is used to reconstruct the corresponding solution, based on some additional quantities that are
computed together with the PAF values. \\

\noindent We also mention here that in case $v$ is a prime number, theorem \ref{Theorem:Generalized_Vertical_Constraint} and its corollaries are not used in the algorithm,
since they are not applicable. \\

\noindent In terms of the algorithm performance, we have noticed that often it suffices to
generate randomly anywhere between $10$ to $100$ million $A$ and $B$ sequences. 
The distributed sorting phase is quite fast, as it proceeds in breaking down the original files into 
several smaller files. The matching phase is also quite fast. The reconstruction phase is of negligible
cost, computationally. \\

\noindent In terms of the algorithm total running time for all D-optimal matrices 
found in this paper, we used approximately $200$ CPU days to generate candidate 
$A$ and $B$ sequences and approximately an additional $100$ CPU days
to perform the distributed sorting and matching phases.

\section{Open cases for D-optimal SDSs with $v < 200$}

\noindent We summarize in the following table all odd integers $v <
200$ and the corresponding D-optimal SDS parameters, pointing out
which cases remain open and which cases are known to exist. In
constructing this table, we consulted information made available to
us in \cite{Orrick:2010}. We also consulted and updated the table in
\cite{Djokovic:1997}.

\noindent Note that for $v=75$ a non-circulant type D-optimal matrix
of order $2\times v = 150$ has been constructed in
\cite{HolzmannKharaghani:1998}. A circulant type D-optimal matrix of
order $150$ is not currently known.

\noindent The asterisk in the ``Existence'' column indicates that such
an SDS with the corresponding parameters is given in the previous
section of this paper.

\noindent Often there are more than one essentially different sets
of solutions of the Diophantine equation $a^2+b^2=4v-2$, for a
specific value of $v$. For each such set, we have elected to use the
positive values of $a$, $b$, such that $a \leq b$, to construct the
SDS parameters $r,s$, which implies that $r \geq s$, in all the SDS
parameters that we present in our table.

D-optimal SDSs with parameters $(v;r,s;\lambda)$, $\lambda = r+s-(v-1)/2$
$$
\begin{array}{|l|l|l|l|c|}
\hline
v & r & s & \lambda & \mbox{Existence}\\
\hline \hline
3 & 1 & 0 & 0 & \mbox{Yes} \\
5 & 1 & 1 & 0 & \mbox{Yes} \\
7 & 3 & 1 & 1 & \mbox{Yes} \\
9 & 3 & 2 & 1 & \mbox{Yes} \\
13 & 6 & 3 & 3 & \mbox{Yes} \\
13 & 4 & 4 & 2 & \mbox{Yes} \\
15 & 6 & 4 & 3 & \mbox{Yes} \\
19 & 7 & 6 & 4 & \mbox{Yes} \\
21 & 10 & 6 & 6 & \mbox{Yes} \\
23 & 10 & 7 & 6 & \mbox{Yes} \\
25 & 9 & 9 & 6 & \mbox{Yes} \\
27 & 11 & 9 & 7 & \mbox{Yes} \\
31 & 15 & 10 & 10 & \mbox{Yes} \\
33 & 15 & 11 & 10 & \mbox{Yes} \\
33 & 13 & 12 & 9 & \mbox{Yes} \\
37 & 16 & 13 & 11 & \mbox{Yes} \\
41 & 16 & 16 & 12 & \mbox{Yes} \\
43 & 21 & 15 & 15 & \mbox{Yes} \\
43 & 18 & 16 & 13 & \mbox{Yes} \\
45 & 21 & 16 & 15 & \mbox{Yes} \\
49 & 22 & 18 & 16 & \mbox{Yes \cite{Djokovic:1997}} \\
51 & 21 & 20 & 16 & \mbox{Yes \cite{Cohn:1992}} \\
55 & 24 & 21 & 18 & \mbox{Yes \cite{FS:2001}} \\
57 & 28 & 21 & 21 & \mbox{Yes \cite{KKS:1991:Infinite-Series}} \\
59 & 28 & 22 & 21 & \mbox{Yes \cite{FKS:2004}} \\
61 & 25 & 25 & 20 & \mbox{Yes \cite{Djokovic:1991}} \\
63 & 29 & 24 & 22 & \mbox{Yes } \star \\
63 & 27 & 25 & 21 & \mbox{Yes \cite{Djokovic:1991}} \\
69 & 31 & 27 & 24 & \mbox{?} \\
73 & 36 & 28 & 28 & \mbox{Yes \cite{KKS:1991:Infinite-Series}} \\
73 & 31 & 30 & 25 & \mbox{Yes \cite{Djokovic:1997}} \\
75 & 36 & 29 & 28 & \mbox{?} \\
77 & 34 & 31 & 27 & \mbox{?} \\
79 & 37 & 31 & 29 & \mbox{Yes \cite{Djokovic:1997}} \\
85 & 39 & 34 & 31 & \mbox{?} \\
85 & 36 & 36 & 30 & \mbox{Yes \cite{Gysin:Thesis:1997}} \\
87 & 38 & 36 & 31 & \mbox{?} \\
91 & 45 & 36 & 36 & \mbox{Yes \cite{KKS:1991:Infinite-Series}} \\
93 & 45 & 37 & 36 & \mbox{Yes } \star \\
93 & 42 & 38 & 34 & \mbox{Yes \cite{Djokovic:1991}} \\
97 & 46 & 39 & 37 & \mbox{Yes} \cite{Djokovic:1997} \\
99 & 43 & 42 & 36 & \mbox{?} \\
\hline
\end{array}
\quad
\begin{array}{|l|l|l|l|c|}
\hline
v & r & s & \lambda & \mbox{Existence}\\
\hline \hline
103 & 48 & 42 & 39 & \mbox{Yes } \star \\
103 & 46 & 43 & 38 & \mbox{Yes } \star \\
111 & 55 & 45 & 45 & \mbox{?} \\
111 & 51 & 46 & 42 & \mbox{?} \\
113 & 55 & 46 & 45 & \mbox{?} \\
113 & 49 & 49 & 42 & \mbox{Yes \cite{Djokovic:1997,Gysin:1997}} \\
115 & 51 & 49 & 43 & \mbox{?} \\
117 & 56 & 48 & 46 & \mbox{?} \\
121 & 55 & 51 & 46 & \mbox{Yes } \star \\
123 & 58 & 51 & 48 & \mbox{?} \\
129 & 57 & 56 & 49 & \mbox{?} \\
131 & 61 & 55 & 51 & \mbox{Yes } \star \\
133 & 66 & 55 & 55 & \mbox{Yes \cite{KKS:1991:Infinite-Series}} \\
133 & 60 & 57 & 51 & \mbox{?} \\
135 & 66 & 56 & 55 & \mbox{?} \\
139 & 67 & 58 & 56 & \mbox{?} \\
141 & 65 & 60 & 55 & \mbox{?} \\
145 & 69 & 61 & 58 & \mbox{?} \\
145 & 64 & 64 & 56 & \mbox{Yes \cite{Djokovic:1997,Gysin:Thesis:1997,GysinSeberry:1998}} \\
147 & 66 & 64 & 57 & \mbox{?} \\
153 & 72 & 65 & 61 & \mbox{?} \\
153 & 70 & 66 & 60 & \mbox{?} \\
157 & 78 & 66 & 66 & \mbox{Yes, \cite{Gysin:Thesis:1997,GysinSeberry:1998}} \\
159 & 78 & 67 & 66 & \mbox{?} \\
163 & 79 & 69 & 67 & \mbox{?} \\
163 & 76 & 70 & 65 & \mbox{?} \\
163 & 73 & 72 & 64 & \mbox{?} \\
167 & 76 & 73 & 66 & \mbox{?} \\
169 & 81 & 72 & 69 & \mbox{?} \\
175 & 81 & 76 & 70 & \mbox{?} \\
177 & 84 & 76 & 72 & \mbox{?} \\
181 & 81 & 81 & 72 & \mbox{Yes \cite{Gysin:Thesis:1997,KSWX:JSPI:1997}} \\
183 & 91 & 78 & 78 & \mbox{Yes \cite{KKS:1991:Infinite-Series}} \\
183 & 83 & 81 & 73 & \mbox{?} \\
185 & 91 & 79 & 78 & \mbox{?} \\
187 & 88 & 81 & 76 & \mbox{?} \\
189 & 92 & 81 & 79 & \mbox{?} \\
189 & 87 & 83 & 76 & \mbox{?} \\
195 & 94 & 84 & 81 & \mbox{?} \\
199 & 93 & 87 & 81 & \mbox{?} \\
 &  &  &  &  \\
 &  &  &  &  \\
\hline
\end{array}
$$

\section{Acknowledgements}
Both authors thank the referees for noticing two minor errors 
and for their constructive comments and suggestions that 
helped improve the paper.
Both authors wish to acknowledge generous support by NSERC.
This work was made possible by the facilities of the Shared 
Hierarchical Academic Research Computing Network (SHARCNET) 
and Compute/Calcul Canada.

\end{document}